\documentclass[12pt]{amsart}
\usepackage{amscd,amssymb}
\usepackage[arrow,matrix]{xy}

\theoremstyle{plain}
\newtheorem{thm}[subsection]{Theorem}
\newtheorem{lem}[subsection]{Lemma}
\newtheorem{prop}[subsection]{Proposition}
\newtheorem{cor}[subsection]{Corollary}

\theoremstyle{definition}
\newtheorem{rk}[subsection]{Remark}
\newtheorem{definition}[subsection]{Definition}
\newtheorem{ex}[subsection]{Example}

\numberwithin{equation}{section}
\newcommand{\OO}{{\mathcal O}}

\newcommand{\G}{{\mathcal G}}
\newcommand{\F}{{\mathcal F}}
\newcommand{\A}{{\mathcal A}}
\newcommand{\B}{{\mathcal B}}
\newcommand{\wM}{\widehat{{M}}}

\newcommand{\wS}{\widehat{{S}}}
\newcommand{\wf}{\widehat{{f}}}

\newcommand{\tf}{\tilde {{f}}}
\newcommand{\tM}{\tilde{{M}}}
\newcommand{\wF}{\widehat{{F}}}
\newcommand{\tL}{\tilde{{L}}}
\newcommand{\wX}{\widehat{X}}

\newcommand{\CC}{{\mathcal C}}
\newcommand{\LL}{{\mathcal L}}

\newcommand{\V}{{\mathcal V}}

\newcommand{\Z}{\mathbb{Z}}
\newcommand{\Q}{\mathbb{Q}}

\newcommand{\C}{\mathbb{C}}

\newcommand{\PP}{\mathbb{P}}

\newcommand{\T}{\mathbb{T}}

\DeclareMathOperator{\Hom}{Hom}

\DeclareMathOperator{\im}{im}

\DeclareMathOperator{\codim}{codim}
\DeclareMathOperator{\supp}{supp}

\begin{document}

\title [Characteristic varieties and constructible sheaves]
{Characteristic varieties and constructible sheaves}

\author[Alexandru Dimca]{Alexandru Dimca }
\address{  Laboratoire J.A. Dieudonn\'e, UMR du CNRS 6621,
                 Universit\'e de Nice Sophia-Antipolis,
                 Parc Valrose,
                 06108 Nice Cedex 02,
                 FRANCE.}
\email
{dimca@math.unice.fr}



\begin{abstract}
We explore the relation between the positive dimensional irreducible components of the
characteristic varieties
of rank one local systems on a smooth surface and the associated (rational or irrational)
pencils. Our study, which may viewed as a continuation of D. Arapura's paper \cite{A}, yields new geometric insight into the translated components relating them to the multiplicities of curves in the associated pencil, in a close analogy to the compact situation treated by A. Beauville  \cite{Beau}.
The new point of view is the key role played by the constructible sheaves naturally arising from local systems.

\end{abstract}

\maketitle

\tableofcontents

\section{Introduction} \label{sec:intro}

Let $M$ be a smooth complex quasi-projective variety. The first characteristic varieties $\V_m(M)$ describe the jumping loci for the dimension of the first twisted cohomology groups $H^1(M,\LL)$, where $\LL$ is a rank one local system on $M$. These characteristic varieties, and their relative position inside the algebraic group
$\T(M)=\Hom(H_1(M),\C^*)$ parametrizing the rank one local system on $M$, play a key role in understanding the fundamental group $\pi_1(M)$, see \cite{DPS1}, \cite{DPS2}.

Since we are interested here only in the first cohomology groups, we can replace $M$ by a
 smooth quasi-projective surface, by taking
a generic linear section and applying a general version of Zariski Theorem, see for instance \cite{D1}, p. 25.
Therefore, to help the reader's intuition, we will assume in this paper that $\dim M=2$, though the results and the proofs hold in any dimension. This is a very interesting setting, as the opening lines of Catanese's Introduction in \cite{Cat2} tell us:

\medskip

"The study of fibrations of a smooth algebraic surface $S$ over a smooth algebraic curve $B$ lies at the heart of the classification theory and of the geometry of algebraic surfaces."

\medskip

The {\it main aim} of this paper is to study the translated components of the
characteristic variety $\V_1(M)$. According to  Arapura's results \cite{A}, such a component $W=W(f,\rho)$ is described by
a pair $(f, \rho)$ where

\noindent (a) $f$ is a surjective morphism $M\to S$, from  the surface $M$ to a smooth curve $S$, having
a connected generic fiber $F$;

\noindent (b)  $\rho \in \T(M)$ is a torsion character such that $W$ is the translate by  $\rho$ of the subtorus $f^*(\T(S))$, i.e. $W=\rho \cdot f^*(\T(S))$.

\medskip

Using the (Logarithmic) Isotropic Subspace Theorem, see Catanese \cite{Cat1} (in the compact case)
and Bauer \cite{Ba} and Catanese \cite{Cat2}, Thm. 2.11 in the quasi-projective case, one can determine in many cases the
various possible morphisms $f$ from certain maximal isotropic subspaces in $H^1(M)$, relative to the cup-product 
$$H^1(M) \times H^1(M) \to H^2(M).$$
A similar approach is provided by the study of the resonance varieties, see \cite{DPS1}.

\medskip

In this paper we assume that the morphism $f$ was already determined and concentrate on the finite order
character $\rho$ above.
Our results can be described briefly as follows.  The characters $\rho$ arising in (b) above for a given map $f$ are parametrized by the Pontrjagin dual
$\hat T(f)=\Hom(T(f), \C^*)$ of a finite group $T(f)$ defined in terms of the topology of the
mapping $f$. This group depends only on the multiple fibers in the pencil associated to $f$, see Theorem \ref{thm5}.
When $\chi(S)<0$, any character in $\hat T(f)$ actually gives rise to
a component, see Proposition \ref{prop2.5}, while for $\chi(S)=0$
(i.e. when $S$ is an elliptic curve, case treated by Beauville \cite{Beau} when $M$ is proper, or when $S=\C^*$) one should discard the trivial character in $\hat T(f)$, see Corollary \ref{cor16}. Moreover, for a generic local system $\LL \in W=W(f,\rho)$,
the dimension $\dim H^1(M,\LL)$ is expressed in terms of the Euler characteristic $\chi(S)$ and
the cardinality of the support of $\rho$, see  Corollary \ref{cor1.15}. 

The case $S=\C^*$ is the most mysterious, and Suciu's example of such a component
for the deleted $B_3$-arrangement given in  \cite{S1}, \cite{S2}
played a key role in our understanding of this question. We consider this component in detail in Examples \ref{exMASTER} and \ref{exfin1}, and the generalization given by the $\A_m$-arrangements, discussed in \cite{C1} and \cite{CDS}, see Examples \ref{exfin2} and \ref{exfin3} below.

Our results are exemplified all along this paper on two types of situations:

\noindent {\bf Case A}, when $M$ is a curve arrangement complement in $\PP^2$, and

\noindent {\bf Case B}, when $M$ is a curve arrangement complement on a normal weighted homogeneous surface singularity, case which includes the Seifert link complements discussed in Eisenbud and Neumann's book  \cite{EN}.

In fact, the reader interested only in Case A may refer to  \cite{D3} for additional information.

\medskip

\medskip
In section 2 we collect some basic facts on regular mappings  $f:M \to S$
and the associated pencils. Lemma \ref{lem1} intends to clarify the key notion of {\it admissible map} used by Arapura in  \cite{A}.

In section 3 we give the main definitions related to characteristic varieties. 
Theorem  \ref{thm1.5} collects some (more or less known) facts on the irreducible components of the characteristic varieties, which are derived by a careful reading of Arapura's paper  \cite{A}.
The key topological property (ii) in Theorem  \ref{thm1.5} was not explicitly stated before (in the proper case,
a related property is used in \cite{Beau}). In Corollary \ref{cor3.2} we give a purely topological proof of this
property.
 

In section 4 we emphasize the key role played in this setting by the constructible sheaves
obtained as direct images of local systems on $M$ under the mapping $f:M \to S$, see for instance
Propositions \ref{prop2.5} and  \ref{prop3.1} and Lemmas \ref{lem3}  and  \ref{lem3.1}.
In particular, for a local system $\LL \in W=W(f,\rho)$,
the dimension $\dim H^1(M,\LL)$ is expressed in terms of the Euler characteristic $\chi(S)$ and
the cardinality of the {\it singular support} $\Sigma(\F)$ of the sheaf $\F=R^0f_*\LL_{\rho}$, see Corollary \ref{cor1.15}.

In the final section we associate to a map  $f:M \to S$ as above a finite abelian group $T(f)$, such that the 
 torsion character $\rho$ is determined by a character $\tilde \rho$ of $T(f)$, see formula  \eqref{k2}.
 We compute this  group $T(f)$ in terms of the multiplicities of the special fibers of  $f$, see Theorem \ref{thm5}. 
 The group $T(f)$ is intimately related to the {\it orbifold fundamental group} $\pi_1^{orb}(f)$ of the map $f$ (and even more so to the corresponding {\it orbifold first homology group} $H_1^{orb}(f)$ of $f$), see Corollary \ref{corB}.
 
 In Theorem \ref{thm16} we show that the character $\tilde \rho \in \hat T(f)$ is trivial if and only if the associated
 constructible sheaf $\F=R^0f_*(\LL_{\rho})$ is a local system on $S$.

We would like to thank Alexander Suciu for interesting and stimulating discussions on the subject of this paper, and for suggesting several improvements
of the presentation.

\section{Generalities on pencils on a surface $M$} \label{sec:two}

Let $\tM$ be a smooth compactification of a complex smooth quasi-projective surface $M$.
Let $f:M \to S$ be a regular mapping, where $S$ is a smooth curve. Then there is a minimal non-empty finite set $A \subset \tM$
such that $f$ has an extention $\tf$ to  $U=\tM \setminus A$ with values in the smooth projective model $\wS$ of $S$. By blowing-up the points in $A$, we pass from $\tM$ to a new compactification $\wM$ of $M$ such that $f$ or $\tf$ is the restriction of a regular morphism $\wf:\wM \to \wS$.

We call any of the morphisms $f$, $\tf$ or $\wf$ above  a {\it pencil of curves}. Such a pencil is {\it rational} if the curve $S$
(or, equivalently, $\wS$) is rational, and it is {\it irrational} otherwise. For any $s \in \wS$, we denote by $\CC_s$ the
corresponding fiber in $\tM$ (obtained by taking the closure of $\tf^{-1}(s)$) or in $\wM$. The corresponding pencil will be denoted sometimes by $\CC=(\CC_s)_{s \in \wS}$.

We recall the following sufficient condition to have a rational pencil. 

\begin{prop} \label{prop2.0}

If the surface $M$ satisfies the condition  $W_1H^1(M,\Q)=0$, where $W$ is the weight filtration of the canonical mixed Hodge structure, and if $f: M \to S$
has a generic connected fiber, then $S$ is a rational curve.
Moreover, the condition  $W_1H^1(M,\Q)=0$ holds if the surface $M$ admits a smooth compactification $\tM$ such that $H^1(\tM,\Q)=0$.
\end{prop} 

To prove this result, we need the following.

\begin{lem} \label{lem1}
Let $X$ and $S$ be  smooth irreducible algebraic varieties, $\dim S=1$ and let $f:X \to S$ be a non-constant
morphism. Then for any compactification $\wf:\wX \to \wS$ of $f$ with $\wX$, $\wS$ smooth, the following are equivalent.

\medskip

\noindent (i) The generic fiber $F$ of $f$ is connected.

\medskip

\noindent (ii) The generic fiber $\wF$ of $\wf$ is connected.

\medskip

\noindent (iii) All the fibers  of $\wf$ are connected.

\medskip

If these equivalent conditions hold, then $f_{\sharp}: \pi _1(X) \to \pi_1(S)$ and $\wf_{\sharp}: \pi _1(\wX) \to \pi_1(\wS)$ are surjective.

\end{lem}

\proof Note that $D=\wX \setminus X$ is a proper subvariety (not necessarily a normal crossing divisor) with finitely many irreducible components $D_m$. For each such component $D_m$, either $\wf(D_m)$ is a point,
or $\wf:D_m \to \wS$ is surjective. In this latter case, it follows that $\dim (\wF \cap D_m) < \dim D_m \leq \dim \wF$. Since $\wF$ is smooth of pure dimension, it follows that $\wF$ is connected if and only if 
$F =\wF \setminus \cup_m ( D_m \cap \wF)$ is connected.
 To show that $(ii)$ implies $(iii)$ it is enough to use the Stein factorization theorem,
  see for instance \cite{H}, p. 280, and the fact that a morphism between two smooth projective curves which is of degree one
 (i.e. generically injective) is in fact an isomorphism.

\medskip
 
 To prove the last claim for $f$, note that there is a Zariski open and dense subset $S' \subset S$ such that
 $f$ induces a locally trivial topological fibration $f:X'=f^{-1}(S') \to S'$ with fiber type $F$. Since $F$ is connected, we get an epimorphism $f_{\sharp}: \pi _1(X') \to \pi_1(S')$. The inclusion of $S'$ into $S$ induces an epimorphism at the level of fundamental groups. Let $j:X' \to X$ be the inclusion. Then
 we have seen that $f \circ j$ induces an epimorphism at the level of fundamental groups. Therefore the same is true for $f$. The proof for $\wf$ is completely similar.

 \endproof

\proof (of Proposition \ref{prop2.0}).
From the surjectivity of $f_{\sharp}$ it follows that  $f^*:H^1(S,\Q) \to H^1(M,\Q)$ is injective. Since $f^*$ preserves the weight filtration $W$, it follows that $W_1H^1(S,\Q)=0$, i.e. $S$ is a rational curve.
 
 If the surface $M$ admits a smooth compactification $j: M\to \tM$ such that $H^1(\tM,\Q)=0$, then 
 $W_1H^1(M,\Q)=j^*H^1(\tM,\Q)=0$ as explained for instance in \cite{D1}, p.243.
 
 \endproof
 
 \begin{ex} \label{ex2.0} The following classes of surfaces $M$ satisfy $W_1H^1(M,\Q)=0$ and will be used as test cases in the sequal.
 
 \medskip

\noindent {\bf Case A}:  complements of plane curve arrangements, i.e. $M=\PP^2 \setminus C$ where $C$ is a plane curve, usually with several irreducible components. One can then take $\tM=\PP^2$. A more general situation is obtained by replacing $\PP^2 $ by any smooth simply-connected surface $\tM$, e.g. $\tM$ a smooth complete intersection in some
projective space.

\medskip

\noindent {\bf Case B}:  complements of curve arrangements on a weighted homogeneous isolated complete intersection $(X,0)$ with $\dim X=2$, whose link  $\Sigma(X,0)$ is a $\Q$-homology sphere.
Such a singularity can by represented by an affine complete intersection surface $X$ defined by some weighted homogeneous equations with respect to some positive integer weights
{\bf w}$=(w_1,...,w_n)$
$$f_1(x)=...=f_{n-2}(x)=0$$ 
in $\C^n$. Moreover $X$ is smooth away from the origin. Let $g \in \C[x_1,...,x_n]$ be another
weighted homogeneous polynomial with respect to the weights
 {\bf w} and set $C=g^{-1}(0)$, $M= X \setminus C$. Since the link of $(X,0)$ is a $\Q$-homology sphere,
one has $H^1(X\setminus 0, \Q)= H^2(X\setminus 0, \Q) =0$. Note that each irreducible component $C_j$ of $C$ is a rational curve, since $C_j^*=C_j \setminus 0$  is
exactly a $\C^*$-orbit of the $\C^*$-action on $\C^n$ associated to the given weights. Moreover, this shows that
$C_j^*$ is smooth. Let $C_1,...,C_r$ be the set of these components.  The Gysin exact sequence
$$ 0= H^1(X\setminus 0   ,\Q) \to H^1(M,\Q) \to \oplus _{j=1,r}H^0(C_j^*,\Q)(-1)\to H^2(X\setminus 0, \Q) =0$$
shows that $H^1(M,\Q)$ is pure of type $(1,1)$, in particular one has $W_1H^1(M,\Q)=0$. Moreover it gives
$b_1(M)=r$.

 \end{ex} 
 
The following explicit description of rational pencils is recalled for the reader's convenience. 

\begin{prop} \label{prop1}
Let $U$ be a smooth surface. 

If $f:U \to \PP^1$ is a morphism, then $L=f^*\OO(1)$ is a line bundle on $U$, generated by the two global
sections  $s_i=f^*(y_i)$, with $y_1,y_2$ a system of homogeneous coordinates on $\PP^1$.

Conversely, if $L$ is a line bundle on $U$, generated by  two global
sections  $s_i$ for $i=1,2$, then there is a morphism  $f:U \to \PP^1$ such that $L=f^*\OO(1)$ and  $s_i=f^*(y_i)$, with $y_1,y_2$ a system of homogeneous coordinates on $\PP^1$.

\end{prop}

\proof

It is well known, see for instance \cite{H}, p. 150, that a morphism $f:U \to \PP^1$ is given
by a line bundle $\LL \in Pic(U)$ and two sections $s_1, s_2 \in \Gamma (U,\LL)$ which do not vanish
both at any point in $U$. In fact $\LL=f^*(\OO(1))$ and $s_i=f^*(y_i)$, with $y_1,y_2$ a system of homogeneous coordinates on $\PP^1$. With this notation, one has $f(x)=[a:b]$ where $[a:b] \in \PP^1$ is such that
$as_2(x)-bs_1(x)=0$.

\endproof

\begin{rk} \label{rk0.3}

Since $U$ is smooth, we have $Pic (U) = C\ell (U)$ and similarly
$Pic (\tM) = C\ell ( \tM)$, see for instance \cite{H}, p. 145. On the other hand, the inclusion
$j:U \to \tM$ induces an isomorphism $j^*:C\ell(\tM) \to C\ell (U)$, as $\codim A =2$, see \cite{H}, p. 133. It follows that $j^*:Pic(\tM) \to Pic (U)$ is also an isomorphism, i.e. any line bundle $L \in Pic(U)$
is the restriction to $U$ of a line bundle $\tL$ on $\tM$. When $\tM=\PP^2$, then $\tL$ has the form 
$\OO(D)$ and the global sections of $L$ are nothing else but the restrictions of global sections of the line bundle $\tL=\OO(D)$, which are the degree $D$ homogeneous polynomials.
In general, the two sections $s_i$ have natural extensions to $\tM$, and we may consider the divisors $C_i:s_i=0$ on $\tM$ and the associated rational pencil $\CC_f: \alpha_1 s_1+\alpha _2 s_2$ of  curves on $\tM$.

\end{rk}

We regard in the sequel the difference $C = \tM \setminus M$ as a (reduced) curve and let $C=\cup _{j=1,r}C_j$, be the decomposition of $C$ into irreducible components. 

\begin{prop} \label{prop2}
Let $B \subset \wS$ be a finite set and denote by $S$ the complement $\wS \setminus B$.
For any surjective morphism $f:M \to S$ and any compactification $\tM$ of $M$ as above,
 any irreducible component $C_j$ of $C$  is in one of the  following cases.

\begin{enumerate}

\item $C_j$ is contained in a curve $\CC_b$ in the pencil $\CC$, corresponding to a point $b \in B$;

\item $C_j$ is strictly contained in a curve $\CC_s$ in the pencil $\CC$, corresponding to a point $ s \in S$;

\item $C_j$ is a horizontal component, i.e. $C_j$ intersects the generic fiber $\CC_t$ of the pencil $\CC$ outside the base locus.

\end{enumerate}

Moreover, if $|B|>1$, then $C_j$ is in the first case above if and only if the homology class $\gamma_j$ of a small loop around $C_j$
satisfies $H_1(f)(\gamma_j) \ne 0$ in $H_1(S,\Z)$.

\end{prop} 

\proof  Let $C_j$ be an irreducible component of $C$.
Then either $\tf(C_j)$ is a point, which leads to the first two cases, or $\tf(C_j)$ is dense in $\wS$, which leads to the last case. The strict inclusion in the second case comes from the surjectivity of $f$.

The last claim is obvious, using the Mayer-Vietoris exact sequence of the covering $\wS=S \cup D$,
where $D$ is the union of small closed discs on $\wS$ centered at the points in $B$. For instance, in the first case, if $\delta _b$ is a small loop at $b$, then one has $H_1(f)(\gamma_j)=m_j\cdot \delta _b$, with $m_j>0$ the multiplicity of the curve
$C_j$ in the divisor $f'^{-1}(b)$ (if the orientations of the loops $\gamma_j$ and $\delta _b$ are properly chosen). See also equation \eqref{eq2} below.

\endproof

\begin{definition} \label{def0}
 In the setting of Proposition \ref{prop2}, we say that the curve arrangement $C$ is {\it minimal} with respect to the surjective mapping $f:M\to S$ if any component $C_j$ of $C$ is of type (1), i.e.
 $C_j$ is contained in a curve $\CC_b$ in the pencil $\CC$, corresponding to a point $b \in B$.
We say that the curve arrangement $C$ is  {\it special} with respect to the surjective mapping $f:M\to S$ if some component $C_j$ of $C$ is of type (2), i.e. $C_j$ is strictly contained in a curve $\CC_s$ in the pencil $\CC$, corresponding to a point $ s \in S$.

\end{definition} 

\begin{rk} \label{rk1}

If $|B|>1$, then the base locus $A$ of the pencil $\CC$ on the surface $\tM$ is just the intersection of any two distinct fibers
$\CC_b \cap \CC_{b'}$ for $b,b'$ distinct points in $B$. Note also that the second case $(2)$ above cannot occur if all
 the fibers $\CC_s$ for $s \in S$ are irreducible. The fibers $\CC_s$ may be non-reduced, i.e. we consider them usually as divisors. Saying that $C_j$ is contained in $\CC_s$ means that $\CC_s=m_j C_j+...$,
with $m_j>0$. 
\end{rk}

 \section{Local systems and characteristic varieties} \label{sec:3}

 \subsection{Local systems on $S$} \label{sec:3.1}
 
 We return to the notation $S=\wS\setminus B$, with $B=\{b_1,...,b_k\}$ a finite set of cardinality
 $|B|=k\ge 0$. Let $g=g(\wS)$ be the genus of the curve $\wS$ and denote by $\delta_1,...,\delta_{2g}$
 the usual $\Z$-basis of the first integral homology group $H_1(\wS)$.
 
 If $\delta_{2g+i}$ denotes an elementary loop based at some base point $b_i \in B$ and turning once around the point $b_i$, then using the usual choices, the first integral homology group of $S$ is given by
\begin{equation} \label{eq2} 
H_1(S)=\Z <\delta _1,...,\delta_{2g+k}>/<\delta _{2g+1}+...+ \delta_{2g+k}>.
\end{equation}
Therefore, for $k>0$, the rank one local systems on $S$ are parametrized by the $(2g+k-1)$-dimensional algebraic torus $\T(S)= \Hom( H_1(S),\C^*)$ given by
\begin{equation} \label{eq3} 
\T(S)=\{ \lambda=(\lambda_1,...,\lambda_ {2g+k})\in (\C^*)^{2g+k}~|~ \lambda _{2g+1}\cdots \lambda_{2g+k}=1\}.
\end{equation}
Here $\lambda_j\in \C^*$ is the monodromy along the loop  $\delta_j$. When $k=0$, one has
\begin{equation} \label{eq1} 
\T(S)= \Hom( H_1(S),\C^*)=\{ \lambda=(\lambda_1,...,\lambda_ {2g})\in (\C^*)^{2g}\}= (\C^*)^{2g}  .
\end{equation}
Note that in both cases $\dim \T(S)=b_1(S)$.
For $ \lambda \in \T(S)$, we denote by $\LL_{ \lambda}$ the corresponding rank one local system on $S$.

The twisted cohomology groups $H^m(S, \LL_{ \lambda})$ are easy to compute. There are two cases.

\noindent Case 1 ($\LL_{ \lambda}=\C_S$). Then we get the usual cohomology groups of $S$, namely
for $k>0$
we have $\dim H^0(S, \LL_{ \lambda})=1$, $\dim H^1(S, \LL_{ \lambda})=2g+k-1$ and $H^m(S, \LL_{ \lambda})=0$
for $m \geq 2$. And for $k=0$ we have $\dim H^0(S, \LL_{ \lambda})= \dim H^2(S, \LL_{ \lambda}) = 1$, $\dim H^1(S, \LL_{ \lambda})=2g$ and $H^m(S, \LL_{ \lambda})=0$
for $m \geq 3$.

\noindent Case 2 ($\LL_{ \lambda}$  nontrivial). This case corresponds to the case when at least one monodromy $\lambda_j$ is not 1. In such a situation one has $2g+k \geq 2$. Then
we have 
\begin{equation} \label{eq1.5} 
\dim H^1(S, \LL_{ \lambda})=2g+k-2=-\chi (S) 
\end{equation}
and $H^m(S, \LL_{ \lambda})=0$ for $m \ne 1$. Indeed, one has obviously $\dim H^0(S, \LL_{ \lambda})=0$ in this case. The vanishing of
$H^2(S, \LL_{ \lambda})$ follows by duality if $k=0$, and since $S$ is homotopically a bouquet of circles
when $k>0$.

\subsection{Local systems on $M$} \label{sec3.2}

The rank one local systems on $M$ are parametrized by the algebraic group
\begin{equation} \label{eq3.5} 
\T(M)= \Hom( H_1(M),\C^*)
\end{equation}
which an extension of the algebraic torus $(\C^*)^{b_1(M)}$ by the finite group $Tors H_1(M)$.
This group can be described explicitely as soon as we know $H_1(M)$.

\noindent {\bf Case A}: complements of plane curve arrangements, i.e. $M=\PP^2 \setminus C$ where $C$ is a plane curve,  with  irreducible components $C_j$ for $j=1,...,r$, $\deg C_j=d_j$.
Let $\gamma_j$ be an elementary loop around the irreducible component $C_j$, for $j=1,...,r.$
Then it is known, see for instance \cite{D1}, p. 102, that
\begin{equation} \label{eq4} 
H_1(M)=\Z <\gamma _1,...,\gamma_r>/<d_1\gamma _1+...+ d_r\gamma_r>
\end{equation}
where $d_j$ is the degree of the component $C_j$.
It follows that the rank one local systems on $M$ are parametrized by the algebraic group
\begin{equation} \label{eq5} 
\T(M)= \Hom( H_1(M),\C^*)=\{ \rho=(\rho_1,...,\rho_r)\in (\C^*)^r ~|~  \rho _1^{d_1} \cdots \rho_r^{d_r}=1\}.
\end{equation}
The connected component $\T^0(M)$ of the unit element $1 \in \T(M)$ is the $(r-1)$-dimensional
torus given by
\begin{equation} \label{eq6} 
\T^0(M)=\{ \rho=(\rho_1,...,\rho_r)\in (\C^*)^r ~|~  \rho _1^{e_1} \cdots \rho_r^{e_r}=1\}
\end{equation}
with $D=G.C.D.(d_1,...,d_r)$ and $e_j=d_j/D$ for $j=1,...,r.$

\begin{rk} \label{rk0}
If $d_1=1$, then $\{\gamma_2,..., \gamma_r\}$ is a basis for $H_1(M)$ and the torus $\T(M)$ can be identified to
$(\C^*)^{r-1}$ under the projection $\rho \mapsto (\rho_2,...,\rho_r)$.

\end{rk}

\noindent {\bf Case B:} complements of curve arrangements on a weighted homogeneous isolated complete intersection $(X,0)$ with $\dim X=2$, whose link  $\Sigma(X,0)$ is a $\Z$-homology sphere.
Using the notation from Example \ref{ex2.0}, we see that $H_1(M)=  \Z <\gamma _1,...,\gamma_r>$,
where $\gamma_j$ is an elementary loop around the irreducible component $C_j$, for $j=1,...,r.$
It follows that
\begin{equation} \label{eq6.5} 
\T(M)= \Hom( H_1(M),\C^*)= (\C^*)^r .
\end{equation}

The computation of the twisted cohomology groups $H^m(M, \LL_{ \rho})$ is one of the major problems.
A simple situation is described in the following.

\begin{ex} \label{ex3.0} 
 
 If $\LL_{ \rho}=\C_M$ and we are in one of the two cases above, this computation can be done as follows. 

 \noindent {\bf Case A}.    The result depends on the local singularities of the plane curve
$C$. In fact $\dim H^0(M, \C)=1$, $\dim H^1(M, \C)=r-1$ and $H^m(M, \C)=0$
for $m \geq 3$. To determine the remaining Betti number $b_2(M)=\dim H^2(M, \C)$ is the same as determining
the Euler characteristic $\chi (M)=3-\chi (C)$ and this can be done, e.g. by using the
formula for $\chi(C)$ given  in  \cite{D1}, p. 162. 

\noindent {\bf Case B}.  Here one has $\dim H^0(M, \C)=1$, $\dim H^1(M, \C)=r$, $\dim H^2(M, \C)=r-1$ and $H^m(M, \C)=0$
for $m \geq 3$. To see this, note that $M$ is an affine open subset in $X$ (which yields $H^m(M, \C)=0$
for $m \geq 3$), and there is a $\C^*$-action on $M$ with finite isotropy groups (which yields $\chi (M)=0$).

\end{ex}

\medskip

To study these cohomology groups $H^m(M, \LL_{ \rho})$  in general, one idea is to study the characteristic varieties
\begin{equation} \label{eq7} 
\V_m(M)=\{ \rho \in \T(M) ~|~  \dim H^1(M, \LL_{ \rho}) \ge m \}.
\end{equation}

\subsection{Arapura's results} \label{sec:3.3}

We recall here some of the main results from \cite{A}, applied to the rank one local systems on $M$, with some additions from  \cite{L1}, \cite{DPS1} and some new consequences.

\begin{thm} \label{thm1.5}

 Let $W$ be an irreducible component of $\V_1(M)$ and assume that  $d_W :=\dim W \ge 1$.
 Then there is a surjective morphism $f_W:M \to S_W$ onto a smooth curve $S_W$, with a connected generic fiber $F(f_W)$,
and a torsion character $\rho_W \in \T(M)$  such that  
$$W= \rho_W \otimes  f_W^*(\T(S_W)).$$
More precisely, the following hold.

\medskip

\noindent (i) $S_W=\wS _W \setminus B_W$, with  $B_W$ a  finite set  satisfying $d_W=2g(\wS _W)+k_W-1=-\chi(S_W)+1$ if
$k_W:=|B_W|>0$. If $B_W = \emptyset$, then $d_W=2g(\wS _W)= -\chi(S_W)+2   $.

\medskip

\noindent (ii)  For any local system $\LL \in W$, the restriction $\LL | F(f_W)$ of $\LL$ to the generic fiber of $f_W$ is trivial, i.e. $\LL | F(f_W)=\C_{ F(f_W)}$.


\medskip

\noindent (iii) If $N_W$ is the order of the character $ \rho_W$, then there is a commutative diagram
$$\xymatrix{
M'  \ar[r]^p \ar[d]^{f'_W} & M \ar[d]^{f_W} \\
S'_W  \ar[r]^q & S_W
}$$
where $p$ is a unramified $N_W$-cyclic Galois covering , $q$ is possibly ramified  $N_W$-cyclic Galois covering, $f'_W$ is $\mu _{N_W}$-equivariant in the obvious sense, has a generic fiber $F(f'_W)$ isomorphic to the  generic fiber $F(f_W)$ of  $f_W$,
and $p^*\rho_W $ is trivial. Here $\mu _{N_W}$ denotes the cyclic group of the $N_W$-th roots of unity.

\medskip

\noindent (iv) If $1 \in W$ and $\LL  \in W$, then $\dim H^1(M,\LL) \ge -\chi(S_W)$ and equality holds with finitely many exceptions.

\medskip

\noindent (v) If $1 \in W$, then $d_W \ge 2$ for $k_W>0$ and $d_W \ge 4$ for $k_W=0$.

\medskip

\noindent (vi) If $1 \notin W$ and either $d_W =2$ and $S_W$ is not an elliptic curve,
or $d_W >2$, then the subtorus $W'= f_W^*(\T(S_W))$ is another irreducible component of $\V_1(M)$.

\end{thm}

\proof

The first claim is just Thm. 1.6 in \cite{A}, section V.

\medskip

To prove the claim (i) note that $d_W=\dim W=\dim \T(S_W)=b_1(S_W)$.

For the claim (vi), consider now the situation $1\notin W$. Note that $\chi(S_W) \geq 0$ if and only if either $g(\wS _W)=0$
and $k_W \leq 2$ or $g(\wS _W)=1$ and $k_W=0$. The first possibility (which contains the trivial cases
$S_W=\PP^1$, $S_W=\C$ and the interesting case $\C^*$) is excluded, since then
$\dim W =\dim \T(S_W) \leq 1$. The second case corresponds to $S_W$ being an elliptic curve $E$, and can be excluded as above if we assume $d_W>2$. The corresponding translated components in this case (assuming $M$ proper) are described in
\cite{Beau}. A uniform treatment of the translated components in the only two
interesting cases $\C^*$ and $E$ is given below in Corollary \ref{cor16}.

With the exception of these special cases, it follows that $\chi(S_W)<0$ and  $W'= f_W^*(\T(S_W))$ is another irreducible component of $\V_1(M)$ by Prop. 1.7 in \cite{A}.

\medskip

Now we prove the claim (ii). Since  $W= \rho_W \otimes  f_W^*(\T(S_W))$, it is enough to prove the claim for 
 $\LL_1=\LL_{ \rho_W }$. Since $\rho_W$ is a torsion character, it follows that
  $\LL_1$ is a unitary local system.  Let $\overline M$
 be a good compactification of $M$ obtained by adding the normal crossing divisor $D$ to $M$. Let $(L_1, \nabla _1)$ be the integrable flat connection on $M$ corresponding to the local system
 $\LL_1$ and let $(\overline L_1, \overline \nabla _1)$ be the Deligne extension of the connection
 $(L_1, \nabla _1)$ to $\overline M$ with residues having the real parts in $[0,1)$. Then there is 
 a Hodge to de Rham spectral sequence
 $$E_1^{p,q}=H^q(\overline M, \Omega_{\overline M} ^p(logD)\otimes \overline L_1) \implies H^{p+q}(M,\LL_1).$$
 Since by hypothesis $H^1(M,\LL_1) \ne 0$, it follows that either $E_2^{1,0} \ne 0$, or $E_2^{0,1} \ne 0$.
 In the first case, we are exactly in the situation of Prop. 1.3 in section V, \cite{A} and our claim (ii) is proved in the final part of the proof. Just note that on the last line of this proof, one should replace ``which forces $\psi |F$ to be trivial''
by  ``which forces $\psi | (F\cap X)$ to be trivial''. (This is due to the fact that $F$ in  \cite{A}
denotes the compactification of our affine fiber $F= F(f_W)$, and $X$ in  \cite{A} corresponds to our $M$.)
If we are in the latter case, then one can show that $\LL_1^{-1}$ leads to the first case, exactly as in the second part of the proof of  Prop. 1.4 in section V, \cite{A}. Since the claim (ii) for  $\LL_1$ is
equivalent to the claim (ii) for $\LL_1^{-1}$, we are done.
The property (ii) corresponds to Prop. 1.2 in Beauville's paper \cite{Beau}. As noted there, it is the same thing to ask triviality for the restriction to one generic fiber of $f_W$ or to all generic fibers of $f_W$.
See Corollary \ref{cor3.2} below for a direct topological proof of the property (ii).

\medskip

The claim (iii) is just the ``untwisting'' part of the proof of Thm. 1.6  in  \cite{A}. The existence of the diagram is explained there via the Stein factorization for $f_W \circ p$. However, the fact that the morphism $q$ has degree $N_W$ depends on the previous claim (ii), and this key point is not mentioned in  \cite{A}.

\medskip

The proof of the claim (iv) is more technical. Using the Projection Formula
\begin{equation} \label{eqp1} 
p_*(\C_{M'})\otimes \LL \simeq p_*(p^*(\LL))
\end{equation}
for $\LL \in W$, see for instance  \cite{D2}, p.42 and then the Leray Spectral Sequence for $p$,
see for instance  \cite{D2}, p. 33, one gets an isomorphism of $\mu _{N_W}$-representations
\begin{equation} \label{eqp2} 
H^1(M',p^*\LL)=H^1(M,p_*(\C_{M'})\otimes \LL ).
\end{equation}
Following the argument in the proof of Thm. 1.6  in  \cite{A}, we get the following
$$\dim H^1(M,\LL) \ge -\chi (S_W).$$
The only point which deserves some attention is the fact that $S_W$ and $\wS _W$ do not admit finite triangulations as claimed  in  \cite{A}, since they are not compact. However, we can replace them by finite
simplicial complexes without changing the homotopy type, e.g. $S_W$ can be replaced by the compact Riemann surface with boundary obtained from $\PP^1$ by deleting small open discs centered at the points in $B_W$.

The fact that there are only finitely many local systems $\LL \in W$ such that 
$$\dim H^1(M,\LL)> -\chi (S_W)    $$ 
follows
by an argument similar to the end of the proof of Prop. 1.7 in  \cite{A}, section V. For a different approach and a generalization to translated components, see  Corollaries \ref{cor1.15}  and \ref{cor6} below.

 \medskip

Finally, the claim (v) follows directly from (iv).

\endproof

\begin{rk} \label{Rk}
Conversely, if $f:M \to S$ is a morphism with a generic connected fiber and with $\chi(S)<0$, then
$W_f=f^*(\T(S))$ is an irreducible component in $\V_1(M)$ such that $1 \in W_f$ and $\dim W_f \ge 2$, see
 \cite{A}, Section V, Prop. 1.7. Some basic situations of this general construction of irreducible components $W_f$ are the following.
 
 \medskip
 
 \noindent {\bf Case A.}
 
 \medskip
 
\noindent (i) {\it The local components}, see for instance  \cite{S1}, subsection (2.3) 
in the case of line arrangements. The case of curve arrangements in $\PP^2$ runs as follows. Let $p \in \PP^2$ be a point such that there is a degree $d_p$ and an integer $k_p >2$ such that

\begin{enumerate}

\item the set $A_p=\{j ~~|~~p \in C_j \text{ and } \deg C_j=d_p\}$ has cardinality $k_p$;

\item $\dim <f_j~~|~~j \in A_p>=2,$ with $f_j=0$ being an equation for $C_j$.

\end{enumerate}
If $\{P,Q\}$ is a basis of this 2-dimensional vector space, then the associated pencil induces a map
$$f_p:M \to S_p$$
where $S_p$ is obtained from $\PP^1$ by deleting the $k_p$ points corresponding to the curves
$C_j$, for $j \in A_p$. In this way, the point $p$ produces an irreducible component in $\V_1(M)$, namely
$$W_p= f^*_p(\T(S_p))$$
of dimension $k_p-1$, and which is called local because it depends only on the chosen point $p$. Note that in the case of line arrangements $p$ can be chosen to be any point of multiplicity at least $3$.

\medskip

\noindent (ii) {\it The  components associated to neighborly partitions}, see  \cite{LY}, corresponds
exactly to pencils associated to the line arrangement, as remarked in  \cite{FY}, see the proof of Theorem 2.4

\medskip

\noindent {\bf Case B.} Let $(X,0) \subset (\C^n,0)$ be as in Example \ref{ex2.0}.
Let $g_1$ and $g_2$ be two weighted homogeneous polynomials of degree $d$ with respect to the weights {\bf w}
such that 
$$X \cap \{g_1=0\} \cap \{g_2=0\}=0.$$
Define $g:X \setminus 0 \to \PP^1$ by $x \mapsto (g_1(x):g_2(x))$. Note that $g$ is constant on the corresponding $\C^*$-orbits. Assume that the generic fiber of $g$ is connected, i.e. it coincides to an orbit. Let $B \subset \PP^1$ be a finite subset such that $k=|B|>2$ and $\CC_b=g^{-1}(b)$ is connected
for any $b \in B$. Then if we set $S=\PP^1 \setminus B$, $\CC =\cup_{b \in B}\CC_b$ and $M=(X \setminus 0) \setminus \CC$,
we have
$H_1(M)=\Z^k$, $\T(M)=(\C^*)^{k}$,  $\T(S)=(\C^*)^{k-1}$ and the subtorus $W=g^*(\T(S))$
is a $(k-1)$-dimensional irreducible component of $\V_1(M)$.

\end{rk}

All these points in Case A. are illustrated by the following beautiful example.

\begin{ex} \label{exMASTER} This is a key example discovered by A. Suciu, see Example 4.1 in  \cite{S1}
and Example 10.6 in  \cite{S2}. Consider the line arrangement in $\PP^2$ given by the equation
$$xyz(x-y)(x-z)(y-z)(x-y-z)(x-y+z)=0.$$

We number the lines of the associated affine arrangement in $\C^2$ (obtained by setting $z=1$) as follows: $L_1: x=0$, $L_2:x-1=0$, $L_3: y=0$, $L_4: y-1=0$, $L_5:x-y-1=0$, $L_6:x-y=0$ and $L_7:x-y+1=0$, see the pictures in  Example 4.1 in  \cite{S1} and Example 10.6 in  \cite{S2}. We consider also the line at infinity
$L_8:z=0$.
As stated in  Example 4.1 in  \cite{S1}, there are

\medskip

\noindent (i) Seven local components: six of dimension 2, corresponding to the triple points, and one of dimension 3, for the quadruple point.

\medskip

\noindent (ii) Five components of dimension 2, passing through 1, coming from the following neighborly
partitions (of braid subarrangements): $(15|26|38)$, $(28|36|45)$, $(14|23|68)$, $(16|27|48)$ and $(18|37|46)$. For instance, the pencil corresponding to the first partition is given by
$P=L_1L_5=x(x-y-z)$ and $Q=L_2L_6=(x-z)(x-y)$. Note that $L_3L_8=yz=Q-P$, is a decomposable fiber in this pencil.

\medskip

\noindent (iii) Finally, there is a 1-dimensional component $W$ in $\V_1(M)$ with 
$$\rho_W=(1,-1,-1,1,1,-1,1,-1) \in \T(M) \subset (\C^*)^8$$
and $f_W:M \to \C^*$ given by
$$f_W(x:y:z)=\frac{x(y-z)(x-y-z)^2}{(x-z)y(x-y+z)^2}$$
or, in affine coordinates
$$f_W(x,y)=\frac{x(y-1)(x-y-1)^2}{(x-1)y(x-y+1)^2}.$$
Then $W \subset \V_1(M)$ and $W \cap  \V_2(M)$ consists of two characters, $\rho_W$ above and
$$\rho _W'=(-1,1,1,-1,1,-1,1,-1).$$
Note that this component $W$ is a translated coordinate component. This is related to the fact that the associated pencil is special. For more on this arrangement see Example \ref{exfin1}.

\end{ex}

\section{Translated components  and constructible sheaves} 

We need the following version of the {\it projection formula}, which is used very often, e.g. \cite{A},
\cite{L1}, but for which I was not able to find a reference.

\begin{lem} \label{lem2}
For any local system $\LL_1$ on $M$ and any local system $\LL_2$ on $S$, one has
$$(Rf_*\LL_1)\otimes \LL_2=Rf_*(\LL_1 \otimes f^{-1}\LL_2).$$
\end{lem}

\proof

To prove this Lemma, we start with the usual projection formula, i.e. with the above notation
\begin{equation} \label{pf1}
(Rf_!\LL_1)\otimes \LL_2=Rf_!(\LL_1 \otimes f^{-1}\LL_2)
\end{equation}
see Thm. 2.3.29, p.42 in \cite{D2}. Let $Z$ be a connected smooth complex algebraic variety of dimension $m$.
Then the dualizing sheaf $\omega _Z$ is just $\C_Z[2m]$ and $D_Z\LL=\LL^{\vee}[2m]$ for any local system $\LL$ on $Z$, see Example 3.3.8, p.69 in \cite{D2}. Note also that for two bounded constructible complexes
$\A^*$ and $\B^*$ in $D_c^b(Z,\C)$ we have the isomorphisms
\begin{equation} \label{dual1}
D_Z\A^*\otimes \B^*=RHom(\A^*,\omega _Z)\otimes \B^*=RHom(\A^*,\omega _Z\otimes \B^*)=
\end{equation}
$$=RHom(\A^*, \B^*)[2m].$$
 It follows that
\begin{equation} \label{dual2}
D_Z(\A^*\otimes \B^*)=RHom(\A^* \otimes \B^*,\omega _Z)=RHom(\A^*, RHom(\B^*, \omega _Z))=
\end{equation}
$$=D_Z\A^* \otimes D_Z \B^*[-2m].$$
For the second isomorphism here we refer to Prop. 10.23, p.175 in \cite{Bo}.
Apply now the duality functor $D_S$ to the projection formula \eqref{pf1}. In the left hand side we get
$D_S((Rf_!\LL_1)\otimes \LL_2)=D_S(Rf_!\LL_1)\otimes D_S(\LL_2)[-2]=Rf_*(D_M\LL_1)\otimes D_S(\LL_2)[-2]=
Rf_*(\LL_1^{\vee}  )\otimes \LL_2  ^{\vee}[4].$
Except the isomorphisms explained above we have used here the isomorphism $D_SRf_!=Rf_*D_M$, see Cor. 4.1.17,
p.90 in \cite{D2}. Similarly, the in the right hand side we get
$D_SRf_!(\LL_1 \otimes f^{-1}\LL_2)=Rf_*D_M(\LL_1 \otimes f^{-1}\LL_2)=Rf_*(\LL_1^{\vee} \otimes (f^{-1}\LL_2 ) ^{\vee})[4].$ Since $(f^{-1}\LL_2 ) ^{\vee}=f^{-1}(\LL_2  ^{\vee})$ and since any local system is the dual of its own dual, the proof is completed.

\endproof

Note that $\F=R^0f_*(\LL_1 )$ is in general no longer a local system on $S$, but a {\it constructible sheaf}.
By definition, it exists a minimal finite set $\Sigma=\Sigma(\F) \subset S$, called the {\it singular support } of $\F$,  such that $\F | (S \setminus \Sigma)$ is a local system, see \cite{D2}, p. 87. The main properties
of this sheaf are given in the following result.

\begin{lem} \label{lem3}
Let $\LL_1$ be a rank one local system on $M$, $F$ the generic fiber of $f:M \to S$ and set 
$\F=R^0f_*(\LL_1 )$. Then either

\medskip

\noindent (i) the restriction $\LL_1|F$ is trivial, $\F | (S \setminus \Sigma)$ is a rank one local system
and $\F_s =0$ if and only if  $s \in \Sigma$, or

\medskip

\noindent (ii) the restriction $\LL_1|F$ is non-trivial and $\F=0$.

\end{lem}

\proof

Consider first the case (i).
If $S' \subset S$ is a Zariski open subset such that the restriction $f':M' \to S'$
with $M'=f^{-1}(S')$, is a topologically locally trivial fibration, it follows that $\F|S'$ is a rank 1 local system.
Indeed, for $s \in S'$ we have
$$\F_s=\lim_{s\in D}\F(D)=\lim_{s\in D}H^0(f^{-1}(D),\LL_1)=\C.$$
Here the limit is taken over all the sufficiently small open discs $D$ in $S$ centered at $s$, and the last
equality comes from the fact that the inclusion $F_s=f^{-1}(s) \to f^{-1}(D)$ is a homotopy equivalence and
$\LL_1|F_s=\C_{F_s}$ (recall that $F_s$ is connected, and hence $f^{-1}(D)$ is connected as well). 
In particular $\Sigma \subset S \setminus S'$, and hence
$\Sigma =\emptyset$ if $f:M \to S$ is a locally trivial fibration. 
The above
argument shows also that $\F_s =0$ if and only if  $s \in \Sigma$.

In the case (ii), 
assume that $\F_s \ne 0$ for some $s \in S$. Then there is a small open discs $D$ in $S$ centered at $s$
such that $H^0(f^{-1}(D),\LL_1)\ne 0$. This implies that the restriction $\LL_1|f^{-1}(D)$ is trivial,
and hence $\LL_1|F$ is trivial as well, a contradiction.

\endproof

We have the following key result.

\begin{prop} \label{prop2.5}
Let $f:M \to S$ be a surjective morphism with a generic connected fiber $F$ from
the surface $M$ onto the curve $S$.
Then for any local system $\LL_1$ on $M$ and any local system $\LL_2$ on $S$, one has the following exact sequence
$$0 \to H^1(S,R^0f_*(\LL_1 )\otimes \LL_2) \to H^1(M, \LL_1 \otimes f^{-1}\LL_2) \to H^0(S,R^1f_*(\LL_1 )\otimes \LL_2) .$$
The last morphism is surjective in any of the following situations:

(i) $S$ is affine;

(ii) $\LL_1|F$ is non-trivial;

(iii) $\LL_1|F$ is trivial and $\LL_2$ is generic, i.e. it is different from a finite set of local systems
depending on $f$ and $\LL_1$.

\end{prop}

\proof  We use the Leray spectral sequence
$$E_2^{p,q}=H^p(S,R^qf_*(\LL_1 \otimes f^{-1}\LL_2))$$
converging to $H^{p+q}(M,\LL_1 \otimes f^{-1}\LL_2)$.  By Lemma
\ref{lem2} we have 
$$R^qf_*(\LL_1 \otimes f^{-1}\LL_2)=R^qf_*(\LL_1 )\otimes \LL_2.$$
In particular, the above spectral sequence yields the following exact sequence
$$
0 \to H^1(S,R^0f_*(\LL_1 )\otimes \LL_2) \to H^1(M, \LL_1 \otimes f^{-1}\LL_2) \to K_2^{0,1} \to 0
$$
where $K_2^{0,1}$ is the kernel of the differential $E_2^{0,1} \to E_2^{2,0}$.

When $S$ is affine, this spectral sequence degenerates at $E_2$ since
$E_2^{p,q}=0$ for $p \notin \{0,1\}$ by Artin Theorem, see Thm. 4.1.26, p. 95 in \cite{D2}, and this proves the claim (i).

In the case (ii) one has $E_2^{2,0}= H^2(S,R^0f_*(\LL_1 )\otimes \LL_2)=0$ since $\F=R^0f_*(\LL_1 )=0$.

For the case (iii), we use the exact sequence of cohomology with compact supports
$$0= H^1( \Sigma,\F \otimes \LL_2) \to H^2_c(U,\F \otimes \LL_2) \to     H^2(S,\F \otimes \LL_2)\to H^2( \Sigma,\F \otimes \LL_2)=0$$
where $U=S \setminus \Sigma$ (note that $S$ can be assumed to be compact, since otherwise we are in the affine case (i)), see for instance \cite{D2}, p.46.
Now $\F|U=\LL_0$ is a rank one local system, and we can use duality to get
$$H^2_c(U,\LL_0 \otimes \LL_2)=H^0(U,\LL_0^{\vee} \otimes \LL_2^{\vee} )^{\vee}.$$
These cohomology groups are clearly trivial for $\LL_2|U \ne \LL_0^{-1}$. Since the restriction
$\LL_2|U$ determines the local system  $\LL_2$, this means that there is at most one local system  $\LL_2$ for which
$E_2^{2,0}\ne 0$.

\endproof

To continue we need the following.

 \begin{lem} \label{lem3.1}
The constructible sheaf $\G=R^1f_*(\C_M)$ has no section with finite support.
 
\end{lem}

\proof

This proof is given in D. Arapura \cite{A}, Proposition 1.7, but we repeat it here for the reader's convenience,
and for clarifying some points in Arapura's proof. Let $D$ be a small disc in $S$ centered at a bifurcation point $b \in S$,
let $D^*=D \setminus \{p\}$ and choose a point $q \in D^*$. Set $X_D=f^{-1}(D)$, $M_D^*=f^{-1}(D^*)$ and $M_q=f^{-1}(q)$.
The claim is equivalent to showing that the morphism
$$i_q^*: H^1(M_D,\C) \to H^1(M_q,\C)$$
induced by the inclusion $i_q:M_q \to M_D$ is injective.
Indeed, one has natural identifications $\G_b=H^1(M_D,\C)$ and $\G_q=H^1(M_q,\C)$ and $i_q^*$ corresponds to the corresponding restriction morphism $\G_b \to \G_q$.
The open inclusion $j_b:M_D^* \to M_D$ induces clearly a surjective morphism $H_1(M_D^*) \to H_1(M_D)$, and hence an injective morphism $j_b^*:H^1(M_D,\C) \to H^1(M_D^*,\C)$.

Now, if the disc $D$ was chosen small enough, the restriction of $f$ over $D^*$ is a locally trivial fibration
with fiber type $M_q$ and hence we get the following exact sequence (which is dual to an exact sequence similar to
\eqref{es1}).
\begin{equation} \label{es7}
0 \to H^1(D^*,\C)  \stackrel{f^*} \longrightarrow      H^1(M_D^*,\C)   \stackrel{\iota_q^*} \longrightarrow      H^1(M_q,\C) 
\end{equation}
where $\iota_q:M_q \to M_D^*$ is the inclusion. It follows that $i_q^*: H^1(M_D,\C) \to H^1(M_q,\C)$ is injective if and only if $I=\im (j_b^*) \cap \im \{H^1(D^*,\C)  \stackrel{f^*} \longrightarrow      H^1(M_D^*,\C)\}=0$.
Since $f:M \to S$ is surjective, it follows that $H=f^{-1}(b)$ is a hypersurface in $M$. Let $p$ be a smooth point on the associated reduced hypersurface. It follows that there is an analytic curve germ $\phi :(\C,0) \to (M,p)$ such that $f(\phi(t))$
has some order $d\ge 1$, where $d$ is the multiplicity of $H$ at $p$. Note that in D. Arapura's proof \cite{A},
the multiplicity $d$ is suppose to be 1, which is not always the case.

Let $\sigma \in I$. Since $\sigma \in \im \{H^1(D^*,\C)  \stackrel{f^*} \longrightarrow      H^1(M_D^*,\C)\}$,
it follows that there is a $\beta \in \Hom (H_1(D^*), \C)=H^1(D^*,\C)$ such that $\sigma=\beta \circ f_*$.
The germ $\phi$ induces a morphism $\phi_*:H_1(D^*) \to H_1(M_D^*)$ such that $f_* \circ \phi_*$ is the multiplication by $d$ on the group $H_1(D^*)=\Z$. It follows that $\sigma \circ \phi_*= d \cdot \beta$.

On the other hand, since $\sigma \in \im (j_b^*)$, there is $\sigma ' \in \Hom(H_1(M_D),\C)$ such that
$\sigma= \sigma ' \circ j_{b*}$. It follows that $\sigma \circ \phi_*=\sigma ' \circ j_{b*}\circ \phi_*$
is trivial, since $j_b\circ \phi$ has an obvious extension $\phi$ from the punctured disc $D^*$ to the disc $D$. In conclusion, $\sigma =0$, and then $I=0$, proving our claim.

\endproof

The result of the above Lemma can be generalized as follows.
 
 \begin{prop} \label{prop3.1}
 Let $f:M \to S$ be a surjective morphism with $\dim S=1$ and a connected generic fiber $F$.
If $\LL$ is a rank one local system on $M$, then the constructible sheaf $\G=R^1f_*(\LL)$ has no section with finite support.
Equivalently,
$$ H^0(S,\G \otimes \LL_2)= 0$$
for all but finitely many local systems $\LL_2 \in \T(S)$.

\end{prop}

\proof

First we check that the two last claims are equivalent.
Locally, the two sheaves $\G$ and $\G \otimes \LL_2$ coincide,
 so they admit in the same time non-zero sections with finite support. If this is the case, then clearly
$$ H^0(S,\G \otimes \LL_2)\ne 0$$
for any local system $\LL_2$.
Suppose now that there are no such sections with finite support. Let $\Sigma':=\Sigma (\G)= \Sigma (\G \otimes \LL_2)$ and note that in this case the restriction 
$$H^0(S,\G \otimes \LL_2) \to H^0(S \setminus \Sigma',\G \otimes \LL_2)$$
is injective.
Since $S \setminus \Sigma'$ is homotopically a bouquet of circles (or a compact curve if $S$ is compact and $\Sigma'=\emptyset$), the last group is  non-zero
exactly when the monodromy of $\LL_2$  along any of the loops forming a basis for the integral homology of $S$ is the inverse of one of the eigenvalues
of the monodromy of the local system $\G|(S \setminus \Sigma')$ along this loop, i.e. for a finite number of local systems $\LL_2$.

With the notation from the proof of Lemma \ref{lem3.1}, we have to prove that the restriction morphism
$$i_q^*: H^1(M_D,\LL) \to  H^1(M_q,\LL|M_q) $$
is injective.

The open inclusion $j_b:M_D^* \to M_D$ induces clearly an epimorphism $\pi_1(M_D^*) \to \pi_1(M_D)$, and hence an injective morphism $j_b^*:H^1(M_D,\LL) \to H^1(M_D^*,\LL)$. This claim follows for instance by using the description of the first twisted cohomology groups $H^1(M,\LL)$ in terms of cross-homomorphisms, see \cite{MM}.

\medskip

\noindent Case 1 (the restriction $\LL |F$ is the trivial local system $\C_F$).

\medskip

To study the local system $\LL'=\LL|M_D^*$, note that it corresponds to a character
$$\rho: \pi_1(M_D^*) \to \C^*.$$
The exact sequence
$$ 1 \to \pi_1(M_q) \to \pi_1(M_D^*) \to \pi_1(D^*) \to 1$$
and the triviality of $\LL|M_q$ (note that $M_q$ is a generic fiber of $f$) imply that $\LL'=f^*(\LL_a)$,
where $\LL_a$ is the rank one local system on $D^*$ with monodromy $a \in \C^*$. For this class of local systems we have a long exact sequence in cohomology
\begin{equation} \label{es9}
 \to H^0(M_q,\C)  \stackrel{h^0-a^{-1}\cdot Id} \longrightarrow      H^0(M_q,\C)   \to 
 H^1(M_D^*,\LL') \stackrel{\iota_q^*} \longrightarrow H^1(M_q,\C) 
\end{equation}
see \cite{D2}, p. 212. Here $h^m$ are the monodromy operators of the fibration $M_q \to M_D^* \to D^*$
and clearly $h^0=Id$ since the fiber $M_q$ is connected.

If $a=1$, then locally at the bifurcation point $b \in S$ we have exactly the same situation as in Lemma \ref{lem3.1}, hence the result is already proven.

If $a \ne 1$, then the morphism $ H^0(M_q,\C)  \stackrel{h^0-a^{-1}\cdot Id} \longrightarrow      H^0(M_q,\C) $
is an isomorphism, which yields an injection $H^1(M_D^*,\LL') \stackrel{\iota_q^*} \longrightarrow H^1(M_q,\C)$.
This gives the result in this case, since the composition of two injections is an injection.

\medskip

\noindent Case 2 ( the restriction $\LL |F$ is a non-trivial local system).

\medskip

In this case $R^0f_*\LL=0$ and the Leray spectral sequence of the fibration $M_q \to M_D^* \to D^*$
yields an isomorphism
$$H^1(M_D^*,\LL) \to H^0(D^*,R^1f_*\LL).$$
Since $H^0(D^*,R^1f_*\LL)$ is just the invariant part of $H^1(M_q,\LL|M_q)$ under the monodromy of the local system $R^1f_*\LL$ on $D^*$, this gives rise to a natural injection
$$H^1(M_D^*,\LL) \stackrel{\iota_q^*} \longrightarrow H^1(M_q,\LL|M_q)$$
which completes the proof in this case as well.

\endproof

The following corollary of the exact sequence in Proposition \ref{prop2.5} and of Proposition \ref{prop3.1}
gives also a new, topological proof for the claim in Theorem \ref{thm1.5} (ii).

\begin{cor} \label{cor3.2}
Let $f:M \to S$ be a surjective morphism with a generic connected fiber $F$ from
the surface $M$ onto the curve $S$ with $b_1(S)>0$.
Then for any local system $\LL_1$ on $M$ such that $\LL_1|F$ is non-trivial, and for any generic local system
$\LL_2 \in \T(S)$, one has $H^1(M, \LL_1 \otimes f^*(\LL_2))=0$.

\end{cor}

As a consequence of Proposition \ref{prop2.5}, we get the following extension of Theorem
\ref{thm1.5},(iv). (This special case corresponds to the case $\LL_{\rho}=\C_M$, when
$R^0f_*(\LL_{\rho})=\C_S$ and hence $\Sigma=\emptyset.$ For an illustration of the general case,
see Example \ref{exfin3}).

\begin{cor} \label{cor1.15}
If $\LL_{\rho}$ is a rank one local system on $M$ such that $\LL_{\rho}|F$ is trivial, then
$$\dim H^1(M, \LL_{\rho} \otimes f^{-1}\LL) \ge -\chi(S)+  | \Sigma (R^0f_*(\LL_{\rho}))|   $$
with equality for all but finitely many local systems $\LL \in \T(S)$. In particular, if $W_{f,\rho}=\rho \otimes f^*(\T(S))$ is a positive dimensional irreducible component of $\V_1(M)$,
then $W_{f,\rho}$ is an irreducible component of $\V_q(M)$,
for any $1 \le q \le q(f,\rho):=-\chi(S)+ | \Sigma (R^0f_*(\LL_{\rho}))|$. Conversely, any positive
dimensional irreducible component of $\V_q(M)$ for $q \ge 1$ is of this type.

\end{cor}

\proof
To estimate $\dim H^1(S,\F \otimes \LL_2)$ we compute
$$\chi(S,\F \otimes \LL_2   )=\dim H^0(S,\F \otimes \LL_2)-\dim H^1(S,\F \otimes \LL_2)= \chi(S \setminus \Sigma) $$
using Thm. 4.1.22, p.93 in \cite{D2}. 
This yields
\begin{equation} \label{dim1}
\dim H^1(S,\F \otimes \LL_2)=\dim H^0(S,\F \otimes \LL_2) -\chi(S)+| \Sigma | \ge -\chi(S).
\end{equation}
In the case $\LL_1=\LL_{\rho}$ such that $\LL_1|F$ is trivial, Proposition \ref{prop2.5} yields
$$ H^1(M, \LL_1 \otimes f^{-1}\LL_2)=H^1(S,R^0f_*(\LL_1 )\otimes \LL_2)$$
for all but finitely many local systems $\LL_2 \in \T(S)$. Similarly, the description of $\Sigma$
given above shows that the group
$H^0(S,\F \otimes \LL_2)$
 is zero unless $\Sigma=\emptyset$ and  $\LL_2= (\F)^{-1}$.

The only thing to explain is the last claim in the case $q>1$. Assume that $W_q$ is a 
positive
dimensional irreducible component of $\V_q(M)$ for $q >1$. Since $\V_q(M) \subset \V_1(M)$, there is
an irreducible component $W$ of $\V_1(M)$ such that $W_q \subset W$. Then the first claim in Corollary
\ref{cor1.15} implies that $W \subset \V_q(M)$, i.e. $W_q=W$.

\endproof

\section{Translated components and multiple fibers} 

Let $W$ be a translated irreducible component  of $\V_1(M)$, i.e. $1 \notin W$. Then, as in Theorem  \ref{thm1.5}, there is a torsion character  $\rho \in \T(M)$ and a surjective morphism $f:M\to S$ with connected generic fiber $F$ such that
\begin{equation} \label{k1}
 W=\rho f^*(\T(S))
\end{equation}
We say in this situation that the component $W$ is associated to the mapping $f$.
In this section we give detailed information on the torsion character $\rho \in \T(M)$ in terms of the geometry of the associated
mapping $f:M \to S$.

\subsection{The general setting} 

Let $F$ be the generic fiber of the mapping $f:M \to S$, i.e. $F$ is the fiber of the topologically locally trivial fibration
$f':M' \to S'$ associated to $f$ as in the previous section. Then, we have an exact sequence
\begin{equation} \label{es1}
H_1(F)   \stackrel{i'_*} \longrightarrow     H_1(M')  \stackrel{f'_*} \longrightarrow      H_1(S') \to 0 
\end{equation}
as well as a sequence 
\begin{equation} \label{es2}
H_1(F)  \stackrel{i_*} \longrightarrow      H_1(M)   \stackrel{f_*} \longrightarrow      H_1(S) \to 0 
\end{equation}
which is not necessarily exact in the middle, i.e. the group
\begin{equation} \label{t}
T(f)=\frac{\ker f_*}{\im i_*}
\end{equation}
is in general non-trivial. Here $i:F\to M$ and $i':F \to M'$ denote the inclusions, and homology is taken with $\Z$-coefficients
if not stated otherwise. 

This group was studied in a compact (proper) setting by Serrano, see \cite{Se},
but no relation to local systems was considered there. On the other hand, this compact situation was also studied by A. Beauville in \cite{Beau}, with essentially the same aims as ours.

The sequence \eqref{es2} induces an obvious  exact sequence
\begin{equation} \label{es3}
0 \to T(f) \to \frac{H_1(M)}{\im i_*} \stackrel{f_*} \longrightarrow      H_1(S) \to 0. 
\end{equation}
Since $ H_1(S)$ is a free $\Z$-module, applying the fuctor $\Hom(-,\C^*)$ to the  exact sequence \eqref{es3}, we get a new exact sequence
\begin{equation} \label{es4}
1 \to \T(S) \to \T(M)_F \to \Hom(T(f),\C^*) \to 1.
\end{equation}
Here $\T(M)_F$ is the subgroup in $\T(M)$ formed by all character $\chi: H_1(M) \to \C^*$ such that
$\chi \circ i_*=0$. This means exactly that the associated local system $\LL _{\chi}$ by restriction to $F$
yields the trivial local system $\C_F$.

The torsion character $\rho \in \T(M)$ which occurs in  \ref{k1} is in this subgroup  $\T(M)_F$, see Theorem  \ref{thm1.5}, (ii).
Moreover, this  character $\rho$ is not unique, but its class 
\begin{equation} \label{k2}
{\tilde \rho} \in \frac{\T(M)_F}{\T(S)} \simeq  \Hom(T(f),\C^*)
\end{equation}
is uniquely determined. From now on, we will regard ${\tilde \rho} \in  \Hom(T(f),\C^*)$. Hence, to understand the possible choices for
 ${\tilde \rho}$, we have to study the group $T(f)$.

\subsection{The computation of the group $T(f)$} 

Let $f:M \to S$ be a surjective morphism with a generic connected fiber $F$ as above. Let $C(f) \subset S$ be a finite, minimal subset such that if we put $S'=S \setminus C(f)$, $M'=f^{-1}(S')$, then the induced mapping $f:M' \to S'$ is a locally trivial fibration. For $c \in C(f)$ we denote by $m_c$ the multiplicity of the divisor $F_c=f^{-1}(c)$.
We have the following result, where the first claim is already in \cite{Beau}, see the remarks after Proposition 1.19, and in Serrano, see \cite{Se}. However, this second author wrongly claims that the isomorphism in (i)
holds for the case (ii) as well. The mistake in \cite{Se} is in the proof of Thm. 1.3, Claim 1, where the relation between the $\gamma_p$'s is incorrect. In the proof below, these 1-cycles $\gamma_p$'s are denoted by
$\delta_c$ and the correct relation is $\Delta=0$.

\begin{thm} \label{thm5} \ \\
\noindent (i) If the curve $S$ is proper, then 
$$T(f) = \left( \oplus_{c \in C(f)}\Z/m_c\Z \right)/(\hat 1, ...,\hat 1).$$

\noindent (ii) If the curve $S$ is not proper, then 
$$T(f) =  \oplus_{c \in C(f)}\Z/m_c\Z .$$

\end{thm}

\proof

The main ingredient to prove this theorem is Lemma 3 in \cite{CKO}, which yields the following exact sequence
\begin{equation} \label{o1}
\pi_1(F) \to \pi_1(M) \to \pi_1^{orb}(f) \to 1.
\end{equation}
Here the {\it orbifold fundamental group} $\pi_1^{orb}(f)$ of the mapping $f$ is the quotient of $\pi_1(S')$ by the normal subgroup generated by the elements $\delta_c^{m_c}$ for $c \in C(f)$, with $\delta_c$ a simple loop going once around the point $c$.
Note that this result is stated in \cite{CKO} under the assumption that the curve $S$ is proper, but the proof given there works for $S$ non-proper as well.

The exact sequence \eqref{o1} yields, by passing to abelianizations, the following exact sequence
\begin{equation} \label{o2}
H_1(F) \to H_1(M) \to H_1^{orb}(f) \to 0.
\end{equation}
We will denote by $f_*^{orb}$ the epimorphism $H_1(M) \to H_1^{orb}(f)$ in the exact sequence above.

Coming back to the notation from subsection (3.1), we get the following presentation for
the {\it orbifold first homology group} $H_1^{orb}(f)$ of the mapping $f$
\begin{equation} \label{o3}
 H_1^{orb}(f)= \Z<\delta_1,..., \delta_{2g+k}; \delta_c \text{ for } c \in C(f)>/<\Delta, m_c\delta_c \text{ for } c \in C(f)>
\end{equation}
where $\Delta= \delta_1+ \cdots + \delta_{2g+k}+ \sum_c \delta_c $. There is a natural surjective morphism
\begin{equation} \label{o4}
 \theta: H_1^{orb}(f) \to H_1(S)
\end{equation}
given by $\delta_i \mapsto \delta_i$ for $i=1,...,2g+k$ and $\delta_c \mapsto 0$ for $c \in C(f)$.
Here we use the presentation for $H_1(S)$ given in the formula \eqref{eq2}. Comparing the exact sequence
\eqref{o2} to the sequence \eqref{es2}, we get an isomorphism
\begin{equation} \label{o5}
 \ker (\theta) \simeq T(f).
\end{equation}
When $S$ is proper we have $k=0$ and the group $\ker (\theta)$ is spanned by the loops $\delta_c$ for $c \in C(f)$,
with the relations $m_c \cdot \delta_c=0$ and $\Delta=\sum_c\delta_c=0$.
This yields the claim (i), since clearly $\Delta$ corresponds to the element $(\hat 1, ...,\hat 1)$.

When $S$ is not proper we have $k>0$ and the group $\ker (\theta)$ is spanned by $\Delta$ and the loops $\delta_c$ for $c \in C(f)$,
with the relations $m_c \cdot \delta_c=0$ and $\Delta=0$. The claim (ii) follows from this description.

\endproof
\begin{cor} \label{corB} There is a non-canonical isomorphism
$$H_1^{orb}(f) \simeq H_1(S) \times T(f).$$
In particular, one has
$$\T^{orb}(f) \simeq \T(S) \times {\hat T(f)}$$
where ${\hat T(f)}=\Hom(T(f),\C^*)$ is the Pontrjagin dual of the finite group $T(f)$ and $\T^{orb}(f)=
\Hom(H_1^{orb}(f),\C^*)$ is the corresponding {\it orbifold character group} of $f$.

\end{cor}

\begin{ex} \label{exB}(The computation of the group $T(f)$ in Case B, the Seifert links)

Let $(X,0)$ be a complex quasi-homogeneous normal surface singularities. Then the surface $X^*=X \setminus \{0\}$ is smooth and it has a $\C^*$-action with finite isotropy groups $\C^*_x$. These isotropy groups can be assumed to be trivial, except for those corresponding to finitely many orbits 
$p_1$,...,$p_s$ in $C(X)=X^*/\C^*$. We set $k_p=|\C^*_p|$ for $p\in P=\{p_1,...,p_s\}$.

The quotient $C(X)$
is a smooth projective curve. 
For any finite subset $B$ in $C(X)$ we get a surjective mapping $f:M \to S$ induced by the quotient map
$f_0:X^* \to C(X)$, where $S=C(X) \setminus B$ and $M=f^{-1}(S)$. 

In addition, the curve $C(X)$ is rational iff the link $L(X)$ of the singularity
$(X,0)$ is a $\Q$-homology sphere (use Cor. (3.7) on p. 53 and Thm (4.21) on p. 66 in \cite{D1}). 
In particular,
if the link $L(X)$ of the singularity
$(X,0)$ is a $\Z$-homology sphere, then  $ H_1(M)=\Z^q$ where $q=|B|$, and a basis is provided by small loops 
$\gamma_b$ around the fiber $F_b=f^{-1}(b)$ for $b \in B$, as explained in subsection \ref{sec3.2}. 

One has $f_*(\gamma_b)=k_b \delta_b$, with $k_b$ the order of the isotropy groups
of points $x$ such that $f(x)=b$, and $\delta_b$ a small loop about $b \in \PP^1$.
The set of critical values of the map $f_0:X^* \to C(X)$ is exactly $P$, and each fiber
$F_p=f_0^{-1}(p)$ is smooth (isomorphic to $\C^*$), but of multiplicity $k_p>1$.
Writing down the map $f_{0*}:H_1(X^*) \to H_1(C(X))$ and using its surjectivity, we get that the integers $k_p$ are pairwise coprime.

Let $(X,0)$ be the germ of an isolated complex surface singularity, 
such that the corresponding link $L_X$ is an integral homology sphere. 
Let $(Y,0)$ be a curve singularity on $(X,0)$. Then using the conic structure
of analytic sets, we see that the local complement $X \setminus Y$, with 
$X$ and $Y$ Milnor representatives of the singularities $(X,0)$ and $(Y,0)$, 
respectively, has the same homotopy type as the link complement 
$M=L_X \setminus L_Y$, where $L_Y$ denotes the link of $Y$.

Moreover, if $(X,0)$ and $(Y,0)$ are quasi-homogeneous singularities 
at the origin of some affine space $\C^N$, with respect to the same 
weights, then the local complement can be globalized, i.e., replaced 
by the smooth quasi-projective variety $X \setminus Y$, where 
$X$ and $Y$ are this time affine varieties representing the germs
$(X,0)$ and $(Y,0)$ respectively.

Using the analytic description of the Seifert link $L=(\Sigma (k_1,...,k_n), S_1 \cup...\cup S_q)$ with $k_j \ge 1$ and $n \ge q \ge 2$  given in \cite{EN}, p. 62 and the above notation, we see that
the link complement $M(L)=\Sigma (k_1,...,k_n) \setminus (S_1 \cup...\cup S_q)$ has the homotopy type of the
surface $M$ obtained from the surface singularity $X$ by deleting the orbits (regular for $k_j=1$ and singular for $k_j>1$) corresponding to the $q$ knots $S_j$, $j=1,...,q$. In other words we have a finite set $B \subset \PP^1$ with $|B|=q$
and a mapping $f:M \to S=\PP^1 \setminus B$. 

Let  $N=k_1 \cdots k_q$, $N_j=N/k_j$ for $1\le j \le q$, $N'=k_{q+1} \cdots k_n$, $N'_j=N'/k_j$ for $q+1\le j \le n$.
We can assume that for $j>q$ one has $k_j=1$ iff $j>q+s$, with $s$ a positive integer. 
The above theorem implies in this case
$$T(f)=\Z/(N'\Z)= \oplus_{q+1\le j \le n} \Z/k_j\Z.$$

\end{ex}

For another way of computing the group $T(f)$ in some cases, we refer to \cite{D3}, Section 6.

\begin{definition} \label{def6}
For a character $\tilde \rho: T(f) \to \C^*$, we
define the support $\supp(\tilde \rho)$ of  $\tilde \rho$ to be the singular set $\Sigma(\F)$ of the constructible sheaf $\F=R^0f_*(\LL_{\rho})$ for some representative $\rho$ of $\tilde \rho$.
\end{definition}

In other words, a critical value $c \in C(f)$ is in $\supp(\tilde \rho)$ if for a small disc $D_c$ centered
at $c$, the restriction of the local system $\LL_{\rho}$ to the associated tube $T(F_c)=f^{-1}(D_c)$ about the fiber $F_c$
is non-trivial.
Since two such representatives $\rho$ differ by a local system in $f^*(\T(S))$, it follows from Lemma \ref{lem2} that this support is correctly defined.

\begin{thm} \label{thm16}
Let $f:M \to S$ be a surjective morphism, with connected generic fiber $F$, and let $\tilde \rho:T(f) \to \C^*$
be a character.
Then the support $\supp(\tilde \rho)$ is empty if and only if the character $\tilde \rho$ is trivial.

\end{thm}

\proof

If the character $\tilde \rho$ is trivial, we can represent it by $\rho=1$ and clearly in this case
$\supp(\tilde \rho)=\Sigma(\C_S)=\emptyset$. 

Conversely, assume now that $\supp(\tilde \rho)=\emptyset$. It follows that for any special value
$c \in C(f)$ and any small tube $T(F_c)$ about the fiber $F_c$, the restriction $\LL_{\rho}|T(F_c)$
is trivial. We know in addition that  $\LL_{\rho}|F$
is trivial for any generic fiber $F$ of $f$.

Let as before $f':M' \to S'$ denote the maximal locally trivial fibration associated to $f$,
and recall that $S'=S \setminus C(f)$. Let $\rho': H_1(M') \to \C^*$ be the composition of the character
$\rho: H_1(M) \to \C^*$ with the morphism $H_1(M') \to H_1(M)$ induced by the inclusion $M' \to M$.
Using the exact sequence \eqref{es1}, it follows that there is a unique character $\alpha':H_1(S') \to \C^*$
such that $\rho'=f'^*(\alpha')$.

Let $c \in C(f)$ be any bifurcation value for $f$ and let $\delta_c$ be the cycle in $H_1(S')$
given by a small loop around $c$. Then, using the fact that $f'$ is a locally trivial fibration
with a connected fiber $F$, it follows that the cycle $\delta_c \in H_1(S')$ has a lifting to a cycle
$\tilde \delta_c \in H_1(M')$ such that $f'_*(\tilde \delta_c)=\delta_c$ and the support of $\tilde \delta_c$
contained in the tube $T(F_c)$. It follows that
$$\rho'(\tilde \delta_c)=1=\alpha'( \delta_c).$$
As a result there is a unique character $\alpha:H_1(S) \to \C^*$, such that $\alpha'$ is the composition of
$\alpha$ with the morphism $H_1(S') \to H_1(S)$ induced by the inclusion $S' \to S$.

Now we replace the representative $\rho$ for $\tilde \rho$ by the character $\rho_1=\rho \cdot f^*(\alpha^{-1})$. It follows that the restriction of $\rho _1$ to $H_1(M')$ is the trivial character.
Using the Mayer-Vietoris sequence to express $H_1(M)$ in terms of the covering $M=M' \cup (\cup_{c\in C(f)}T(F_c)$ we get that the character $\rho_1$ itself is trivial. This clearly implies that the character $\tilde \rho$ is trivial.

\endproof

The following result, based on Corollaries \ref{cor1.15}, \ref{corB}  and Theorem \ref{thm16}, clarifies the case of translated components.

\begin{cor} \label{cor16}
 Let $f:M \to S$ be a surjective morphism, with connected generic fiber $F$.
 
 \noindent(i) If $\chi(S)<0$, then the irreducible components in  $\V_1(M)$ associated to $f$ form a subgroup  in $\T(M)$, isomorphic to the orbifold character group $\T^{orb}(f)$. More precisely, they are given by
 $\hat f^{orb}_*(\T^{orb}(f))$, where the injective morphism $\hat f^{orb}_*:\T^{orb}(f) \to \T(M)$ is the dual of the epimorphism $f_*^{orb}$.

 \noindent(ii) If $\chi(S)=0$, then the irreducible components in  $\V_1(M)$ associated to $f$ are given by
 $ \hat f^{orb}_* (\T^{orb}(f)^*)$, where $\T^{orb}(f)^*$ is obtained from the orbifold character group  $\T^{orb}(f)$
 by deleting the identity connected component.

\end{cor}

The same proof as above yields the following result, to be compared with Theorem \ref{thm1.5}, (iv).

\begin{cor} \label{cor6}
Let $f:M \to S$ be a surjective morphism, with connected generic fiber $F$, such that $\chi(S)\le 0$.
Then, for any character $\tilde \rho: T(f) \to \C^*$, one has
$$\dim H^1(M,\LL_{\rho} \otimes f^*\LL)\ge -\chi(S)+| \supp(\tilde \rho)|           $$
for any local system $\LL \in \T(S)$ and
the above inequality
is an equality for all except finitely many local systems $\LL$.

\end{cor}

For the proofs of the following related two results we refer to \cite{D5}. 

\begin{prop} \label{prop7}
For $f:M \to S$  a surjective morphism, with connected generic fiber $F$, and for a non-trivial element $\tilde \rho$ in the Pontrjagin dual ${\widehat T(f)}$, one has a natural adjunction isomorphism
$$\F=Rj_*j^{-1}\F$$
where $\F=R^0f_*(\LL_{\rho})$ and $j:S \setminus \Sigma(\F) \to S$ is the inclusion. 
In particular, the local system
$j^{-1}\F$ on $S \setminus \Sigma(\F)$ is non trivial.

\end{prop}

\begin{cor} \label{cor8}
With the above notation, if $S$ is a compact curve, then $|\Sigma(\F)| \ne 1.$

\end{cor}

 \begin{ex} \label{exfin1}(The deleted $B_3$-arrangement)
We return to Example \ref{exMASTER} and apply the above discussion to this test case.
The corresponding mapping $f:M \to \C^*$ has $B=\{0, \infty \}$ and $C(f)=\{1\}$. Indeed, with obvious notation, we get the following divisors:
$D_0=L_1+L_4+2L_5$, $D_{\infty}=L_2+L_3+2L_7$ and $D_1=L_6+2L$ where $L:x+y-1=0$ is exactly the line from the $B_3$-arrangement that was
deleted in order to get Suciu's arrangement. Moreover, the associated fibration $f':M' \to S'$ in this case is just the fibration of the $B_3$-arrangement discussed in \cite{FY}, Example 4.6.

The line $L$ is the only new component that has to be deleted, therefore $m''(1)=2$.
Since none of the fibers $D_c$ (in our case there  is just one, for $c=1$) is multiple,  Theorem  \ref{thm5} implies that
$$T(f)=\Z/2\Z.$$
Let $\gamma_i =\gamma(L_i)$. We know that $\rho (\gamma_i)=\pm 1$ and to get the exact values we proceed as follows.
First note that we can choose  $\rho (\gamma_1)= 1$, since the associated torus is
$$f^*(\T(\C^*))=\{(t,t^{-1},t^{-1},t,t^2,1,t^{-2},1)~~ |~~ t \in \C^*\}.$$
(In fact the choice $\rho (\gamma_1)= -1$ produces the character $\rho'_W$ introduced in Example \ref{exMASTER}.)
Next let $\alpha= \sum _{i=1,7}\alpha_i \gamma_i \in H_1(M)$. Then $\alpha \in \ker f_*$ if and only if
\begin{equation} \label{e1}
\alpha_1+\alpha_4+2\alpha_5=\alpha_2+\alpha_3+2\alpha_7.
\end{equation}
In our case, the morphism $\theta:\ker f_* \to \Z/2\Z$ is given by $\alpha \mapsto \alpha_2+\alpha_3 -\alpha_6$
It follows that $\gamma_6 \in \ker f_*$ and $\theta(\gamma_6)=1 \in \Z/2\Z$. It follows that $\rho (\gamma_6)=-1$.

Next $\gamma_1+\gamma_2  \in \ker f_*$ and  $\theta(\gamma_1+\gamma_2    )=1 \in \Z/2\Z$. It follows that
$\rho(\gamma_1) \rho(\gamma_2)=-1$, i.e. $\rho(\gamma_2)=-1$. The reader can continue in this way and get 
the value of $\rho=\rho_W$ given above in  Example \ref{exMASTER}.

\end{ex}

 \begin{ex} \label{exfin2}(A more general example: the  $\A_m$-arrangement) 
Let $\A_m$ be the line arrangement in $\PP^2$ defined by the equation
$$x_1x_2(x_1^m-x_2^m)(x_1^m-x_3^m)(x_2^m-x_3^m)=0.$$
This arrangement is obtained by deleting the line $x_3=0$ from the complex reflection arrangement
associated to the full monomial group $G(3,1,m)$ and was studied in  \cite{C1} and in \cite{CDS}. The deleted $B_3$-arrangement studied above is obtained by taking $m=2$.

Consider the associated pencil
$$(P,Q)=(x_1^m(x_2^m-x_3^m),x_2^m(x_1^m-x_3^m)).$$
Then the set $B$ consists of two points, namely $(0:1)$ and $(1:0)$, and the set $C$ is the singleton
$(1:1)$, see for instance  \cite{FY}, Example 4.6. It follows that $m'(c)=1$, $m''(c)=m$ and hence
via Theorem \ref{thm5} we get
$$T(f)=\frac {\Z}{m\Z}.$$
Using Corollary \ref{cor16}, we expect $(m-1)$  1-dimensional components in $\V_1(M)$, and this is precisely what has been proved in  \cite{C1}, or in  Thm. 5.7 in \cite{CDS}.
There are $r=2+3m$ lines in the arrangement, and to describe these components we use the
coordinates 
$$(z_1,z_2,z_{12:1},...,z_{12:m}, z_{13:1},...,z_{13:m},  z_{23:1},...,z_{23:m})$$
on the torus $(\C^*)^r$ containing $\T(M)$. Here $z_j$ is associated to the line $x_j=0$, for $j=1,2$,
and $z_{ij:k}$ is associated to the line $x_i - w^kx_j$, where $i,j=1,3$, $k=1,...,m$, and
$w=\exp (2\pi {\sqrt -1}/m)$.
All the above 1-dimensional components have the same associated 1-dimensional subtorus
$$\T=f^*(\T(\C^*))=\{(u^m,u^{-m},1,...,1,u^{-1},...,u^{-1},u,...,u) ~~|~~ u \in \C^* \}$$
where $f:M \to \C^*$ is the morphism associated to the pencil $(P,Q)$, and each element
$1$, $u^{-1}$ and $u$ is repeated $m$ times.
Let $\gamma_c$ be an elementary loop about one line $L$  in the fiber $\CC_c$, with multiplicity 1, e.g. $L:x_1-x_2=0$. Similarly, let  $\gamma_b$ be an elementary loop about one line $L'$  in the fiber $\CC_b$, with multiplicity 1, where $b= \infty= (0:1)$, e.g. $L':x_2-x_3=0$.
And let  $\gamma_0$ be an elementary loop about one line $L_0$  in the fiber $\CC_0$, with multiplicity 1, where $0= (1:0)$, e.g. $L_0:x_1-x_3=0$.
One can show easily that

\noindent (i) the classes $[\gamma_c]$ and  $[\gamma_b + \gamma_0  ]$ in the group 
$T(f)$ are independent of the choices made;

\noindent (ii) $[\gamma_c] = - [\gamma_b + \gamma_0  ] $ is a generator of  $T(f).$

It follows that a torsion character $\rho \in \T(M)$ such that $\LL_{\rho}|F=\C_F$ and inducing a nontrivial character 
$\tilde \rho: T(f) \to \C^*$ is given by
$$\rho =(1,1, w^k,...,w^k, w^{-k}, ..., w^{-k}, 1,...,1)$$
for $k=1,...,m-1.$ Here $\tilde \rho ( [\gamma_c])=w^k$ and $\rho$ is normalized by setting the last
$m$ components equal to 1.

\end{ex}

 \begin{ex} \label{exfin3}(A non-linear arrangement)
Consider again the pencil $\CC: (P,Q)=(x_1^m(x_2^m-x_3^m),x_2^m(x_1^m-x_3^m))$ associated above to the $\A_m$-arrangement, for $m \ge 2$.
 We introduce the following new notation: $C=\{ (0:1),(1:0), (1:1)\}$. Let $B \subset \PP^1$ be a finite set such that
$|B|=k\ge 2$ and $B \cap C=\emptyset.$ Consider the curve arrangement in $\PP^2$ obtained by taking the union of the $3m$ lines given by
$$(x_1^m-x_2^m)(x_1^m-x_3^m)(x_2^m-x_3^m)=0$$
with the $k$ fibers $\CC_b$ for $b \in B$. Let $M$ be the corresponding complement and $f:M \to S:=\PP^1 \setminus B$
be the map induced by the pencil $\CC$. Then one has the following.

\noindent (i) $T(f)= \frac{\Z}{m\Z} \oplus \frac{\Z}{m\Z} \oplus   \frac{\Z}{m\Z}.$ Let $e_j$ for $j=1,2,3$ denote
the canonical basis of $T(f)$ as a $ \frac{\Z}{m\Z}$-module.

\noindent (ii) For a character $\tilde \rho: T(f) \to \C^*$, let $W_{\rho}=\LL_{\rho} \otimes f^*(\T(S))$ be the associated  component. Then $\dim W_{\rho}=k-1$ and for a local system $\LL \in W_{\rho}$ one has
$$\dim H^1(M,\LL) \ge k-2+\epsilon (\rho)$$
where equality holds for all but finitely many  $\LL \in W_{\rho}$ and 
$$\epsilon (\rho) = |\{j ~~|~~\tilde \rho (e_j) \ne 1\}| \in \{0,1,2,3\}.$$
Indeed, the set $\{j ~~|~~\tilde \rho (e_j) \ne 1\}$ can be identified with the support $\supp(\tilde \rho)$ and the claim follows from Corollary \ref{cor16} and Corollary \ref{cor6}.
This shows that the various translates $W_{\rho}$ of the subtorus $W'=\T_W= f^*(\T(S))$ have all the same dimension, but they are irreducible components of various characteristic varieties $\V_q(M)$, 
with $q=q(f,\rho)=k-2+\epsilon (\rho)$ as in Corollary \ref{cor1.15}, a fact apparently not noticed before.

\end{ex}

\end{document}